\documentclass{article}
\usepackage{booktabs}
\usepackage{float}
\usepackage{enumitem}
\usepackage{amsmath,amssymb,wrapfig,color,amsfonts}
\usepackage{caption}
\usepackage{graphicx}

\begin{document}
\date{}
\title{Delta Redistribution: a variation of State Redistribution for steady-state embedded boundary computations}
\author{Marsha Berger\footnote{Center for Computational Mathematics, Flatiron Institute, Simons Foundation}  \hspace{1in} Andrew Giuliani $^*$}

\maketitle

\section{Introduction}\label{sec:intro}
State redistribution (SRD) \cite{origSRD} is a recent  method 
to stabilize the evolution of
cut cells on an embedded boundary Cartesian grid. 
Among its nice features is ease
of implementation, high-order extensions, GKS stability for a 1D model
linear advection problem,  and a monotonicity-preserving
variation \cite{stableSRD}. SRD has been used in finite
volume and discontinuous Galerkin codes, and has been implemented in
AMReX \cite{GAetal}
for a variety of physical applications and pdes. 

One drawback of SRD is that if applied to a pde for steady-state
computation, it does not preserve the exact steady-state of the underlying
discretization method. It modifies it in the same way that a splitting
method does -- the final answer  will differ by truncation error. One could use a higher
order SRD to reduce this error, but it won't be zero.
At steady state the base scheme residual is nonzero, and one has to instead
examine  the magnitude of  $\max \, \lvert \mathbf{q}^{n+1}-\mathbf{q}^n\rvert $.

There may be cases where the base scheme has some properties that should be preserved,
e.g. well-balancing,
without the modifications of SRD.
For such cases this note proposes
another approach to solve the small cell problem. It is a variation
of SRD that we call Delta Redistribution (DRD) and  
that converges exactly to the underlying scheme.
DRD is also GKS stable for the 1D model problem (not shown here).
What it loses is the monotonicity property of weighted SRD
mentioned above. 
DRD is similar to flux
redistribution (FRD) \cite{chern:colella,PBCCW95} in its steady-state behavior and conservation
properties, but is more accurate. We point out that the original 
SRD weights are also not monotone but are well-behaved in practice, as are the 
weights used in AMReX.

In the next section we briefly present the DRD
algorithm and its finite-volume implementation, which is almost identical to that of SRD.
This is followed by 2D examples for both smooth and shocked steady-state
flow.  For a smooth flow we show a 
convergence plot so that DRD and SRD can be compared. 
We show several examples of shocked flow, steady and unsteady.
This short note assumes the reader is already familiar with 
finite volume schemes and the small cell problem.

\section{The DRD Algorithm}\label{sec:drd}

Both DRD and SRD start with a {\em provisional} finite volume update that is accurate but
unstable on the cut cells. The update should be as accurate on the cut cells as it is on the full
volume cells in order for redistribution to preserve the accuracy. We describe the first
order accurate version of normal merging in two space 
dimensions
in which each cut cell $C_{i,j}$  with volume
fraction less than the threshold temporarily combines with a neighboring cell $C_{i',j'}$
in its most normal direction away from the embedded
boundary.\footnote{This is the most common case; for examples of other
configurations see the references.} 
This cell is in $(i,j)$'s merging neighborhood, $M_{i,j}$.  In words first,
SRD forms a convex combination
of the provisional states, whereas DRD forms a convex combination of the {\em updates}.
The same weights are used in both SRD and DRD. They can be the original weights 
\cite{origSRD}, or the monotone weights \cite{stableSRD}.

To be precise, write the provisional cut cell update as
\begin{equation}
\hat{U}_{i,j} = U_{i,j}^n - \frac{\Delta t}{V_{i,j}} \sum_{k \in {nbrs}}
\, {\mathbf F_k
\cdot \mathbf n_k} \,  A_k = U_{i,j}^n - D_{i,j}  ,
\end{equation}
where $D_{i,j}$ incorporates the entire update, $V_{i,j}$ is the cut cell volume, and $A_k$
is the face length between cell $(i,j)$ and its neighbor in the update stencil.
In the most common, simplest case there is only one neighbor in the merging neighborhood,
but there may be several until a volume threshold of cells in $M_{i,j}$ is reached.
Typically the threshold is half the full cell volume. The formulas are:
\begin{table}[H]
\centering
\caption{\textbf{Merging Formulas for SRD and DRD\\
                  (original weights)}}
\label{eqn:combo}
\begin{tabular}{c c}
\toprule
method  &  merging formula  \\[.08in]
SRD & $\displaystyle Q^{srd}_{i,j} = \frac{1}{\hat{V}_{i,j}}  \sum_{(r,s)\in M_{i,j}} \;
\frac{V_{r,s}}{N_{r,s}} \, \hat{U}_{r,s}$ \\[.17in]
DRD & $\displaystyle Q^{drd}_{i,j} = \frac{1}{\hat{V}_{i,j}}  \sum_{(r,s)\in M_{i,j}} \; 
\frac{V_{r,s}}{N_{r,s}} \, {D_{r,s}} $\\[.08in]
\bottomrule
\end{tabular}
\end{table}

\noindent Here, $\hat{V}_{i,j} = \sum_{(r,s)\in M_{i,j}} \; \frac{V_{r,s}}{N_{r,s}}$
is the same convex combination as in Table \ref{eqn:combo}.

As is evident, the formulas are the same but applied to either the provisional state
$\hat{U}$ or update $D$.
In these formulas  $N_{r,s}$ = number of merging neighborhoods  that cell ${r,s}$ 
is {\em contained in}.
Call this set $W_{i,j}$, so $N_{ij} = \lvert W_{i,j} \rvert $ is its size.
We call this the overlap count.\footnote{These counts can
be determined in a preprocessing step if the geometry is stationary. Only one of sets M or
W needs to be formed for implementation.}
There is a reciprocal relationship between M and W. 
If my cell is in your
merging neighborhood, then you are counted in my overlap count.
Formally, if $(r,s) \in M_{i,j}$, then $(i,j) \in W_{r,s}$. 

To get to the next time step, the merging neighborhoods are distributed back to the
Cartesian grid in
Table \ref{eqn:back}:
\begin{table}[H]
\centering
\caption{\textbf{Redistribution Formulas}}
\label{eqn:back}
\begin{tabular}{c c}
\toprule
method  &  redistribution formula \\[.08in]
SRD & $\displaystyle U_{i,j}^{n+1} = \frac{1}{N_{i,j}}  \sum_{(r,s)\in W_{i,j}} \,
{Q}^{srd}_{r,s} \qquad\quad $ \\[.17in]
DRD & $\displaystyle U_{i,j}^{n+1} = U_{i,j}^n -  \frac{1}{N_{i,j}}  \sum_{(r,s)\in W_{i,j}} \,
Q^{drd}_{r,s} $\\[.08in]
\bottomrule
\end{tabular}
\end{table}

Both the merging and redistribution formulas show only the original weights.
The monotone weights of \cite{stableSRD} may be used in their place; DRD is
not monotone with either set, although the monotone weights reduce the overshoot.
For the second order version (shown in the examples), we reconstruct a linear function
through weighted centroids, and evaluate it at the centroid of the Cartesian
grid. For details see  \cite{stableSRD}. For
shocked flow, the gradient reconstruction step in both the base and 
redistribution  schemes are limited. The only
difference in the limiting step between SRD and DRD is a negative density and pressure test
that is omitted for DRD. Both SRD and DRD call the same procedure, with an
argument that is either $\hat{U}$ or $D$ on the Cartesian grid, and perform the final
update algebra that either includes $U_{i,j}^n$ or not.

\section{Smooth Flows}
We first show two examples of smooth flows so we can compute convergence rates.
The finite volume scheme uses central differences for gradients on full
cells with regular stencils. On the cut and irregular cells,
we do a first order linear reconstruction of
the gradients using a least squares fit of the neighbors in
a 3-by-3 neighborhood that are in the fluid. No limiters are applied for
these smooth cases. 

\subsection{Steady-state supersonic vortex}
\begin{figure}[h!]
\centering
\includegraphics[height=2.0in]{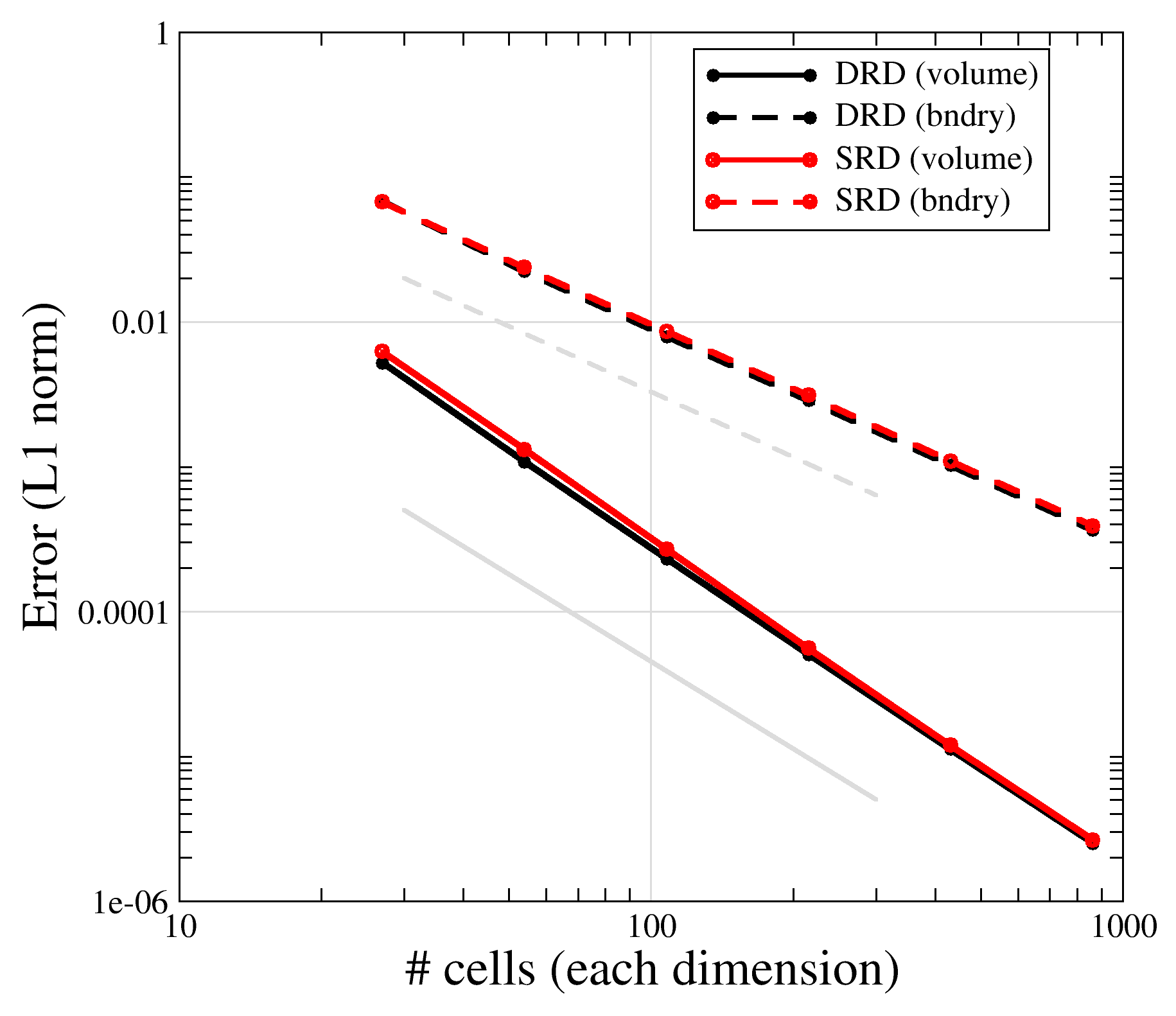}
\hspace*{.5in}
\includegraphics[height=2.05in]{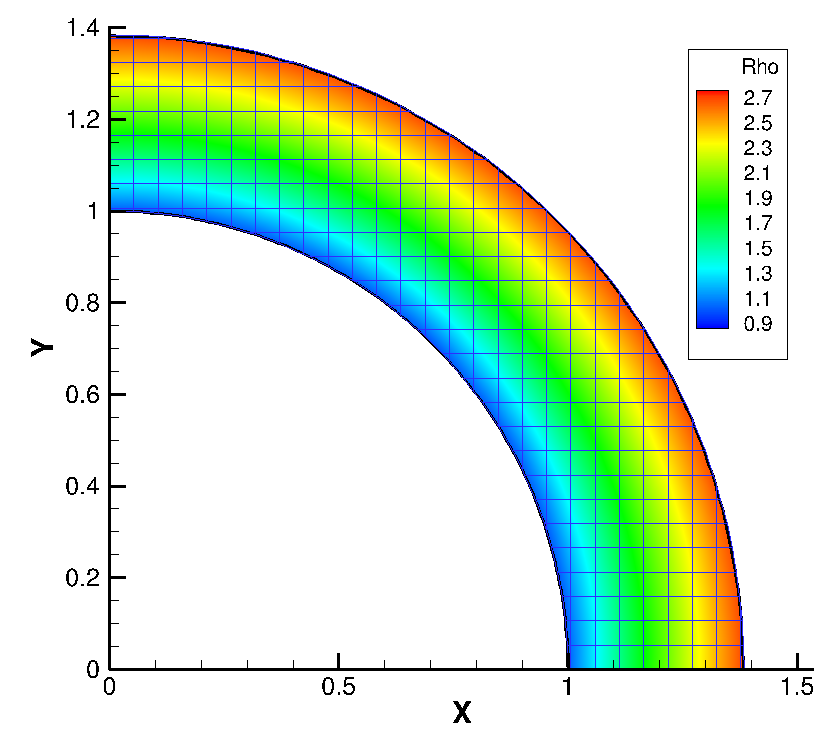}
\caption{Density convergence plot (left) for supersonic vortex test
problem, with density profile shown on the right. 
The convergence rates are second order in L1 for the volume, and 1.5 for
the boundary. The lines showing this are in gray.
\label{fig:ssv}}
\end{figure}

This first example uses the well-known test case of a supersonic vortex in a quarter
circle, originally found in \cite{aftosmis:acc}. We initialize the domain
with the exact solution. We
also use the exact solution at the inflow and outflow boundaries at each stage.
We compare the accuracy of  a linearly exact base scheme, stabilized with
DRD at each timestep, and SRD with linear redistribution. The base
scheme gradient reconstruction is first order accurate and uses 
flow cell neighbors in a 3-by-3
box centered on each cut cell. Regular cells with a regular stencil use
central differences. 
For the redistribution the gradient also uses a 3-by-3 box, extended by
one if necessary if there is a direction without enough cells for a
stable reconstruction. In the steady state case, since DRD converges to the base scheme 
there is no point in
using linear redistribution instead of piecewise constant, which saves
over half the cost.
(For time-dependent simulations linear redistribution is needed to maintain accuracy.)
\footnote{Experimentally, linear and piecewise constant DRD converge to steady state 
at the same rate.} In any case since redistribution is only applied at the
cut cells it is not expensive.
Figure \ref{fig:ssv} shows the convergence rate of the density error in the volume, and
separately the boundary cells.
The L1 norm of the error in the volume is $\sum \lvert (q_{i,j} - exact(x_i,y_j)) \rvert
A_{i,j} $. 
The L1 boundary norm is the density error at cut cell centroids, weighted by the
length of the embedded boundary segment in that cell.
 
\subsection{Linear advection}
DRD is designed for steady-state calculations, but for completeness
we show a time-dependent case of 
linear advection $u_t + a\,u_x + b\,u_y = 0$ along a $39^\circ$ ramp with 
$(a,b) = (\cos(39^\circ),\sin(39^\circ))$. We compare DRD
with and without gradient reconstruction as well as SRD with gradient
for the redistribution.
Since a linear function has a constant update, and DRD is applied to the update,
the piecewise constant
version in the tables above maintains second order accuracy.
However the piecewise linear version has a smaller error constant, which this
test also shows. (The table in the Appendix presents the linear redistribution
step).

\begin{figure}[h!]
\centering
\includegraphics[height=2.1in]{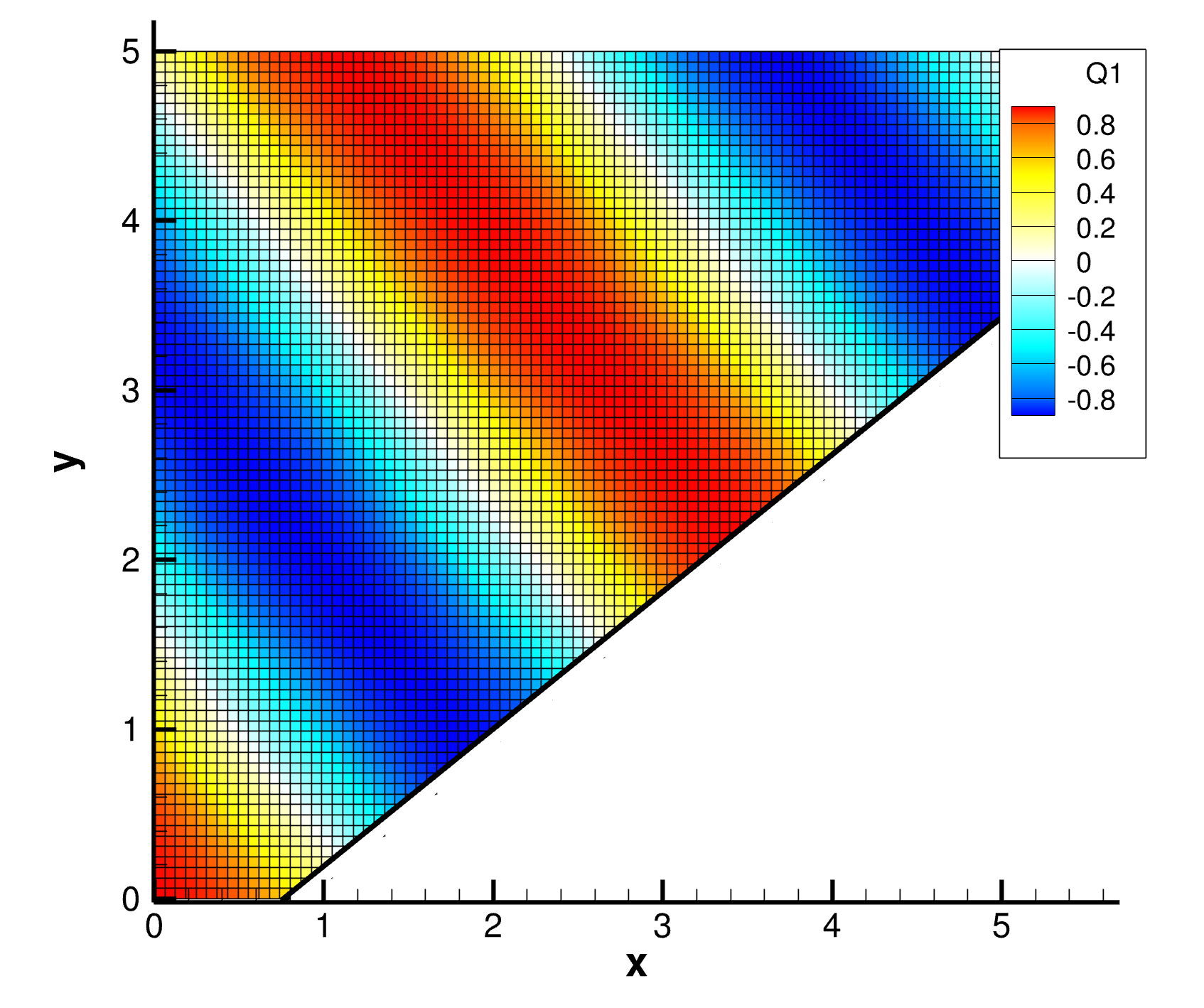}
\hspace*{.5in}
\includegraphics[height=2.1in]{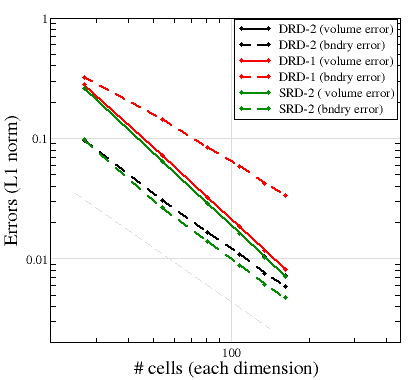}
\caption{Cosine initial conditions for unsteady ramp test (left).
Convergence plot (right)  after integrating to time 2.0.
All three methods converge quadratically in the volume,
and vary between 1.3 and 1.7 along the wall. The dashed gray line shows a
convergence rate of order 1.5.
\label{fig:advfig}}
\end{figure}

The initial conditions are a cosine that travels parallel to the boundary and is
constant normal to the boundary, 
shown in Figure \ref{fig:advfig}(left) on an 81-by-81 mesh.
The computations use a 2 stage Runge-Kutta method to integrate to time 2.0.
Figure \ref{fig:advfig} right shows convergence plots for 
DRD with piecewise constant (DRD-1) and linear
(DRD-2) redistribution, and linear SRD (SRD-2). 
The $L_1$ error in the volume for  DRD-2 and SRD-2 are
indistinguishable, differing in the second or third digits. DRD-1 is 
worse enough that the volume error is also affected, although its
convergence order is still two in the volume. 

\section{Shocked  Flows}
We show three examples of shocked flow, two steady and one moving, to demonstrate DRD's
performance with discontinuities. The limiters applied in each case are
mentioned below.

\subsection{Regular Reflection}
\begin{figure}[b]
\centering
\hspace*{-.2in}
\includegraphics[trim=10 0 0 05,clip,height=1.5in]{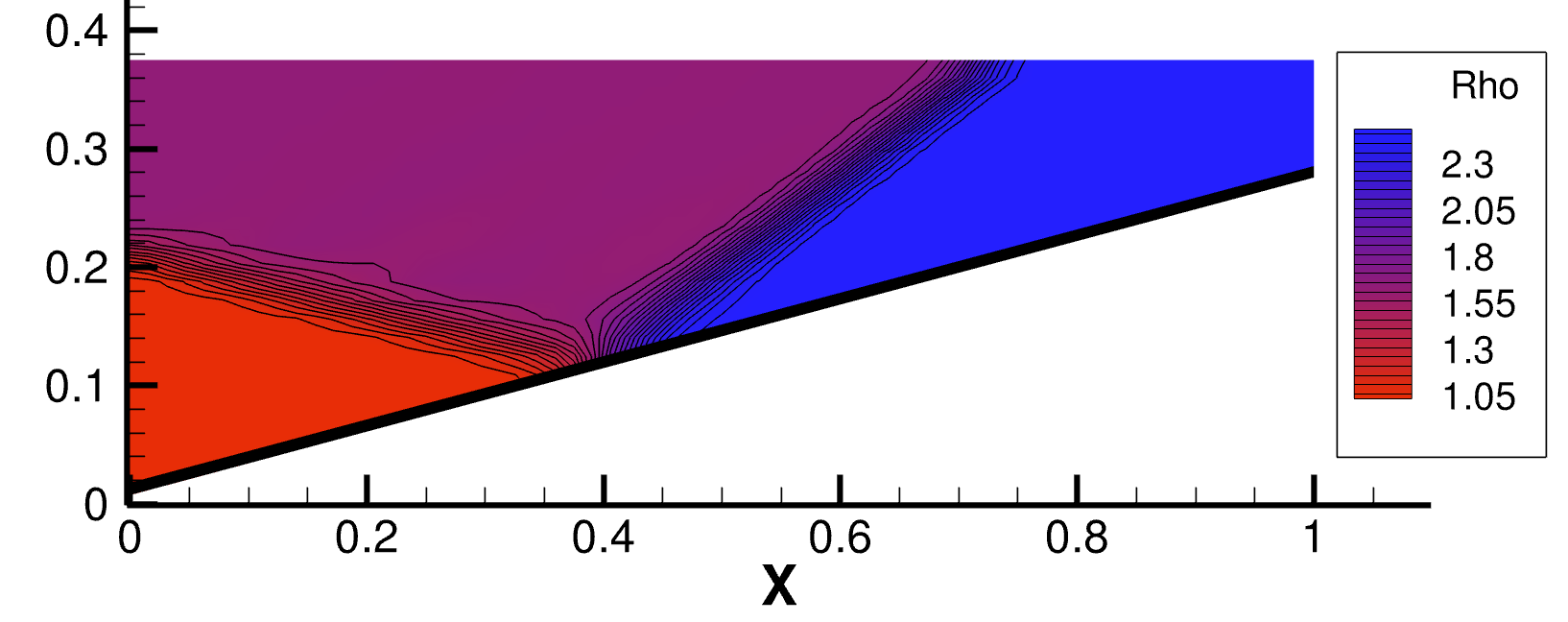}
\hspace*{.2in}
\includegraphics[height=1.5in]{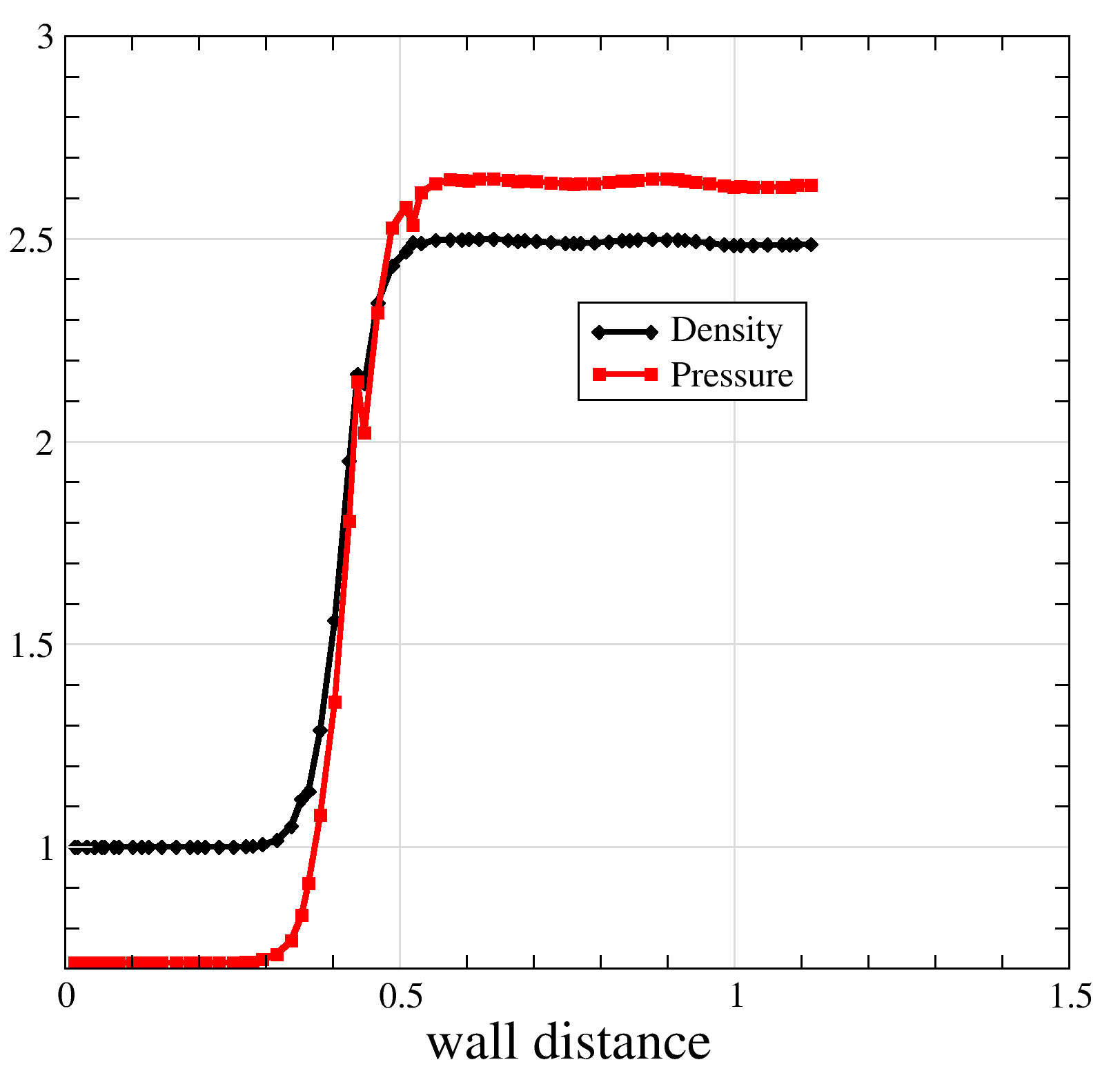}
\caption{Regular reflection at steady-state  computed with DRD. Left is density in the volume, right is
reconstruction along ramp. 
\label{fig:ramp}}
\end{figure}

For a steady shocked flow example figure \ref{fig:ramp} shows a regular reflection on a
$15^\circ$ ramp,
with an inflow Mach number of 2.9. 
The computation is initialized with the exact solution. 
The gradient in the volume uses central differencing with the MC limiter. The cut cells
reconstruct the gradient as above in the smooth case, and limit using Barth-Jespersen.
Limiter chatter begins after about  6 orders of magnitude of convergence (max norm) and 8 orders
(L1 norm). This is enough converged to examine the solution, so 
we do not examine this further  
since the goal is to show stability of DRD for a shocked flow 
problem.  For this steady test we again use DRD without gradients 
for stabilization. The mesh has $\Delta x=\Delta y =1/64$.  The smallest
cut cell in this mesh has a volume fraction of  $0.60\times10^{-3}$. 

\subsection{Steady Bow Shock}
\begin{figure}
\centering
\includegraphics[height=2.2in]{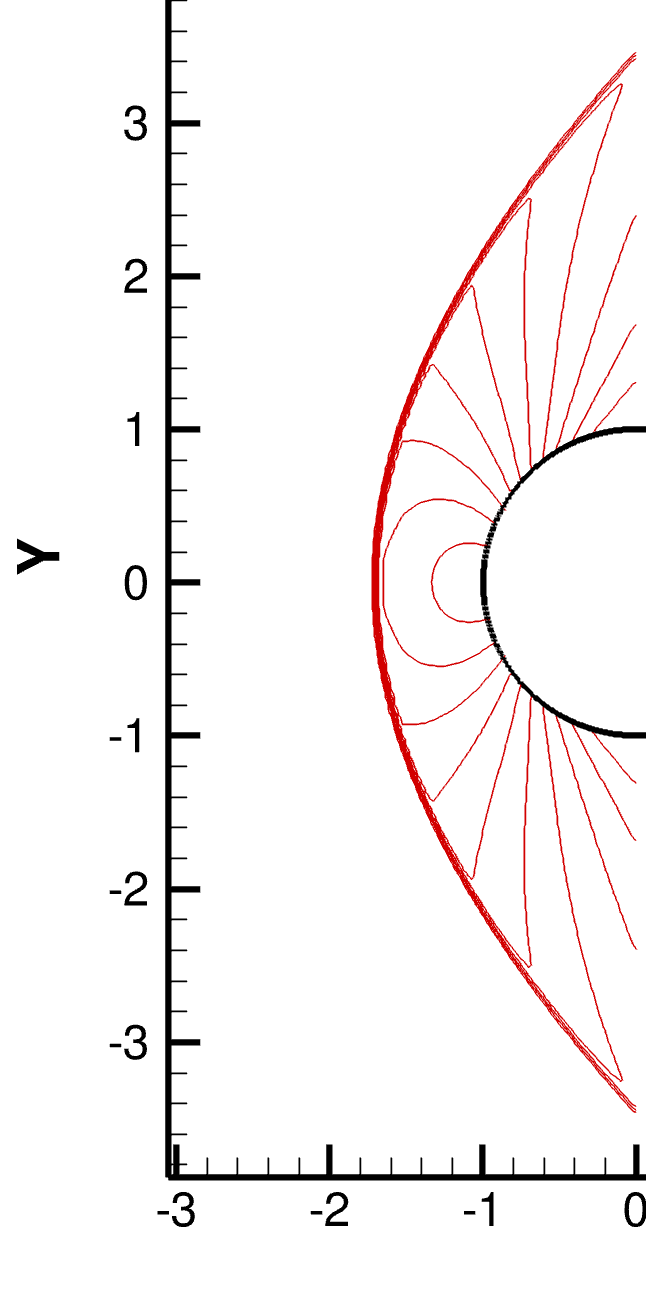}
\hspace*{1.4in}
\includegraphics[height=2.2in]{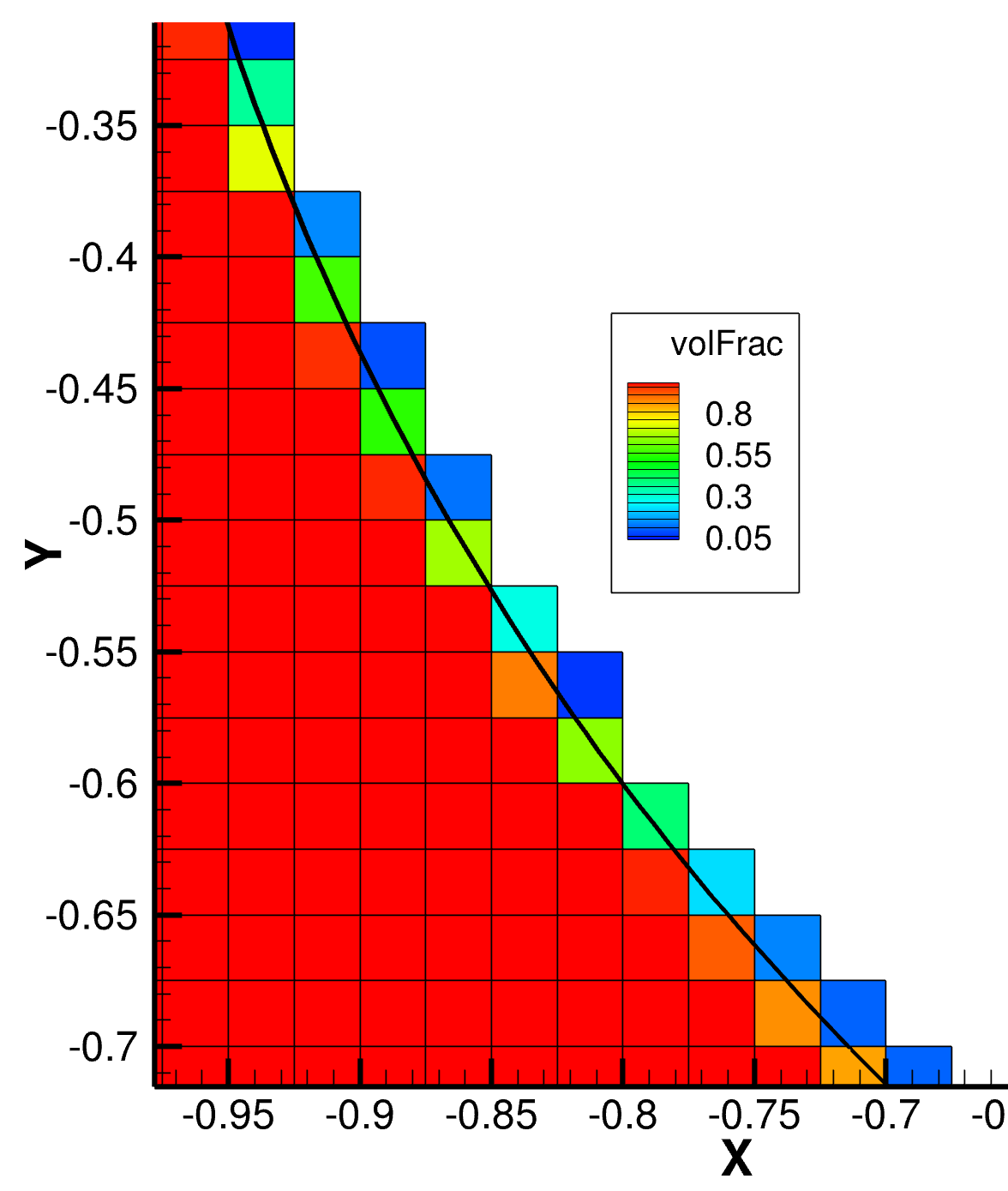}
\caption{Steady blunt body flow computed with DRD. The pressure
contours are smooth. The smallest volume fraction in this grid with
$\Delta x= \Delta y=1/40$ is $0.14\times10^{-2}$. We plot the full Cartesian 
cell even if it is cut to make it easier to visualize.
\label{fig:cylForeBody}}
\end{figure}

The next example shows  Mach 3 flow around a
half cylinder. This example demonstrates that DRD can effectively
compute shocks on varying geometries  as well as 
flat ramps.\footnote{The cut cells in both circular geometries 
treat the wall boundary as a linear segment, the walls are not actually curved.}
We use uniform freestream initial conditions.
The domain and mesh spacing is the same as in  \cite{tanShu}. 
For this example to converge we used the Venkat limiter with K=1 on the
full cells, which includes the shock, and K=10 on the cut cells. Both
values are set to improve  convergence, but in opposite directions: a larger K on the
full cells would not limit the shocks enough.  Using K=1 on the
cut cells does not converge at all.
Figure \ref{fig:cylForeBody} shows the result with  $\Delta x = \Delta y =
1/40$, using the Van Leer Hanel flux function.
The total enthalpy $H = (E+p)/\rho$ should be identically 7 in this
example.
The pre-shock state is exactly this. There is a nonlinear spike across
the shock. The post-shock $H$ is essentially flat
at 6.9925, 6.9966 and 6.9987 on the $h=1/20$, $1/40$ and $1/80$ grids, an error converging
at a rate of  approximately 1.3.  An additional drop over the last few cells next to
the wall brings $H$ to 6.980, 6.993 and 6.996; this wall contribution
converges at roughly first order and so does not degrade the overall rate.
At $h=1/80$ the density and stagnation pressure at the nose are within
0.08\% and 0.13\% of their exact values of 4.307485 and 8.614981,
converging at 1.9 and 1.5 respectively over the finer pair of grids.

\subsection{Unsteady Shocked Flow}
\begin{figure}[h]
\centering
\includegraphics[trim=0 0 5 0,clip,height=2.2in]{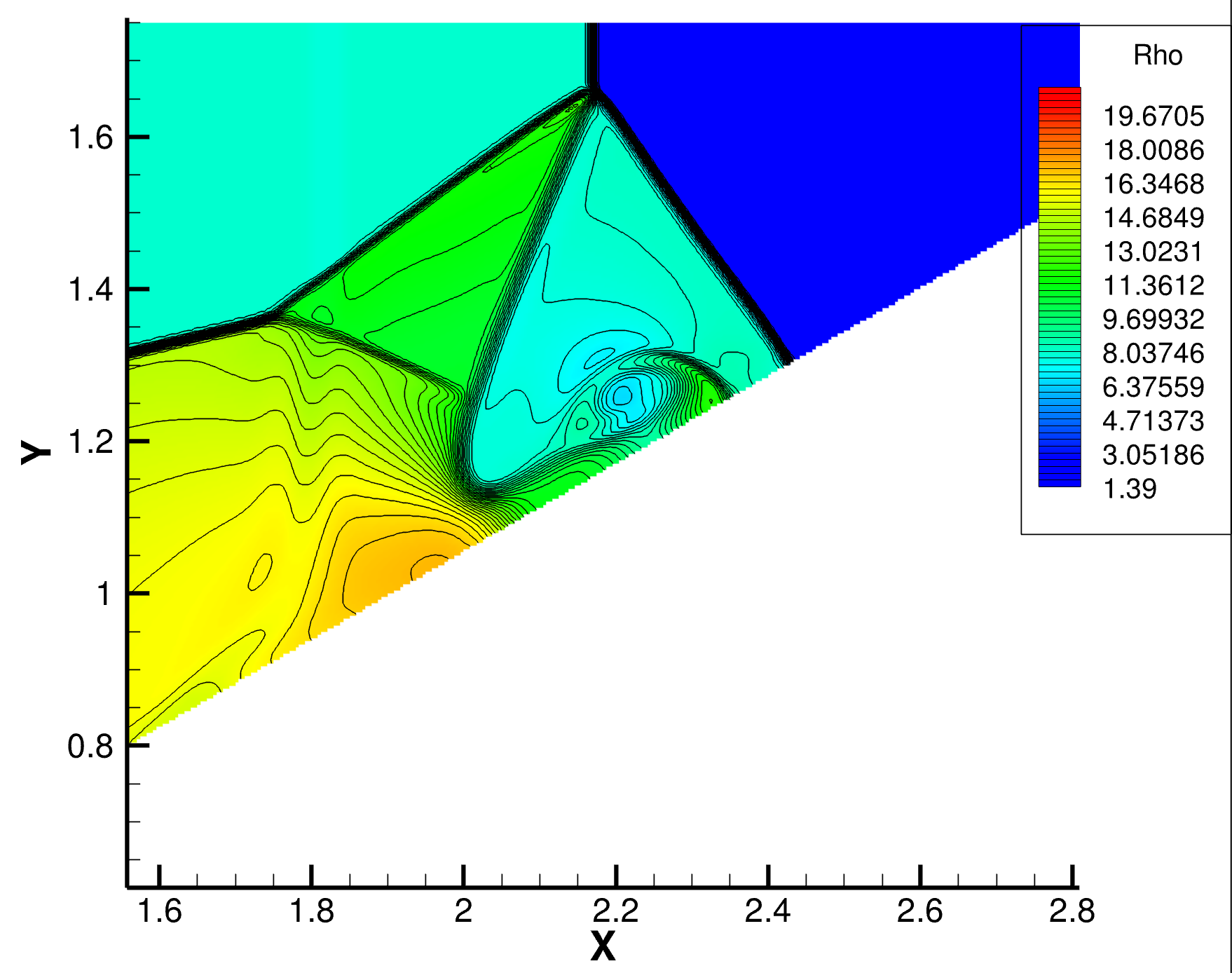}
\hspace*{.1in}
\includegraphics[height=2.2in]{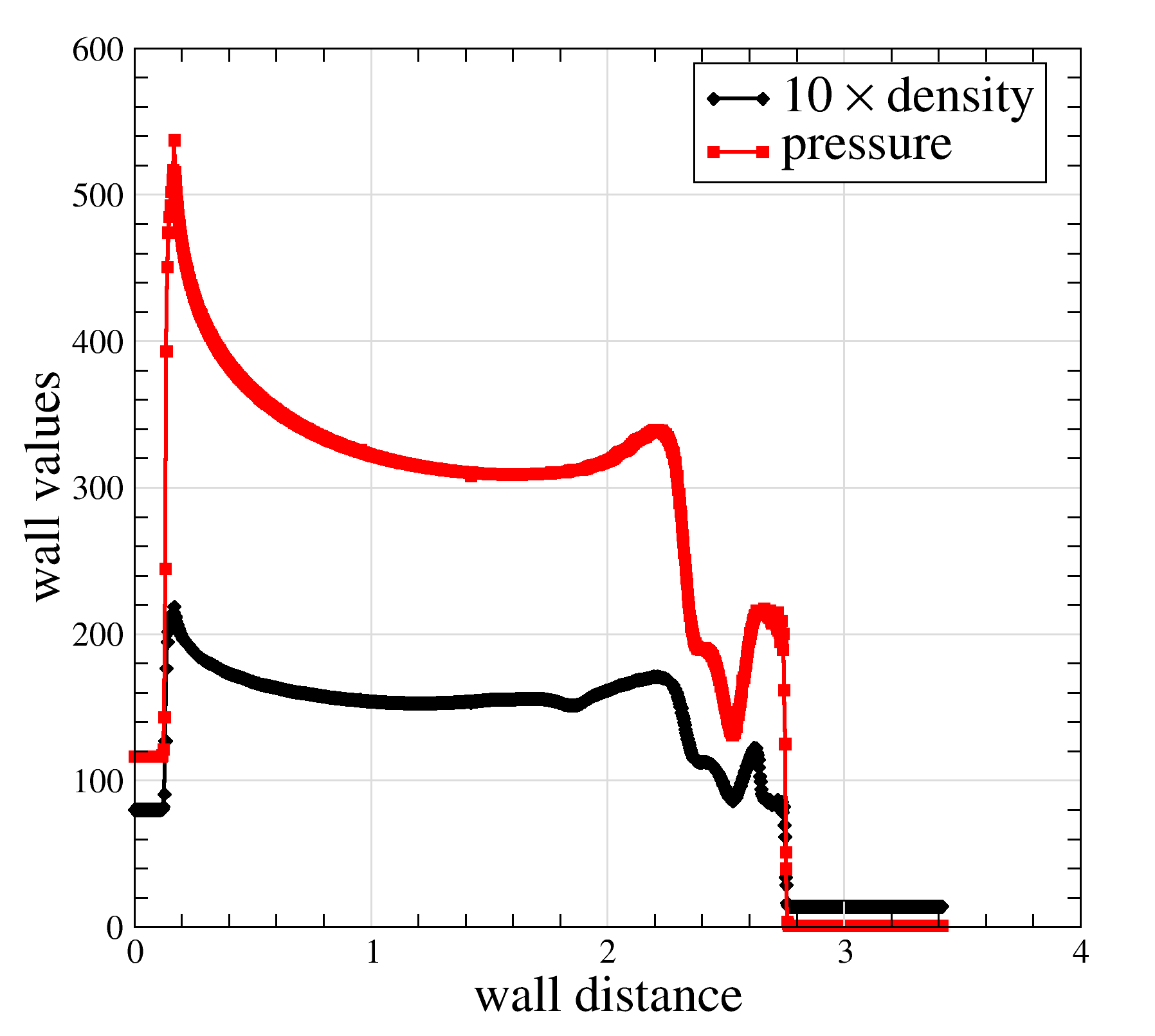}
\caption{Zoom of unsteady Mach 10 shock using DRD (left), and density
and pressure along the embedded wall(right). 
\label{fig:ramp30}}
\end{figure}

We also include here an unsteady example of the Mach 10 shock on a  $30^\circ$ 
ramp.  This is to demonstrate that DRD can successfully handle a moving shock as well.
For this example we use the Colella version of the MUSCL scheme \cite{Colella:Unsplit}, 
modified for cut cells as described in \cite{origSRD}.
Figure \ref{fig:ramp30} (left) shows a zoom of the solution at time 2.0 on a Cartesian
mesh with spacing $h=1/240$. Both the mach stem and contact are seen to be well-preserved.
On the right is the pressure and density (multiplied by 10 so it is
visible)  along the wall, after reconstruction from the cut cell centroids to the
wall centroids.

\section{Conclusion}
This note has presented a method for stabilizing steady flows with cut cells called Delta
Redistribution that is a simple variation of State Redistribution.
Its
purpose is to preserve steady-state solutions of the underlying base
method, where it is cheaper than SRD since no gradients are needed. 
When used with gradients preliminary results show it is effective for 
time-dependent problems as well.
Although it lacks the monotonicity bounds of weighted SRD, it is GKS stable for a linear
advection model problem, as is SRD with the original weights. 
We understand what limiting the reconstruction of a state does to solution
accuracy.  We would like to develop a better understanding of what 
limiting the update does.

\section{Acknowledgements/Use of Generative AI}
This research did not receive any specific grant from funding agencies
in the public, commercial, or not-for-profit sectors. We acknowledge the
use of Claude, which checked our Fortran results in Matlab, suggested
exact solutions for  the steady shock problems, and proofread
the paper.

\section{Appendix}
The following steps compute gradients for $Q_{i,j}$ in
merging neighborhood $M_{i,j}$ for redistribution:
\begin{itemize}[leftmargin=*]
\item
Compute weighted centroids analogous to the weighted
volume, repeated here:
\begin{align*}
\hat{V}_{i,j} &= \sum_{(r,s)\in M_{i,j}} \; \frac{V_{r,s}}{N_{r,s}}\\
\hat{\mathbf{X}}_{i,j} &= \frac{1}{\hat{V}_{i,j}}
   \sum_{(r,s)\in M_{i,j}} \; \frac{V_{r,s} {\mathbf{X}}_{r,s}}{N_{r,s}}
\end{align*}
where ${\mathbf{X}}_{i,j} = (x_{i,j},y_{i,j})$ is the Cartesian cell centroid.
\item
Compute merged  values ${Q}$ as usual for DRD or SRD using Table
\ref{eqn:combo}.
\item
Compute a linear function $q_{i,j}$ 
through $Q_{i,j}$ using least squares
with slopes $s_{x_{i,j}},s_{y_{i,j}}$ fit to the weighted neighborhood centroids:
\begin{equation}
{q}_{i,j}(x,y) = {Q}_{i,j} + s_{x_{i,j}}\,(x-\hat{x}_{i,j}) + s_{y_{i,j}}\, (y-\hat{y}_{i,j})
\end{equation}
For stability, the stencil should include  at least one  centroid 
no closer than $h/2$ in each direction.
\item
Evaluate ${q}$ at the {\em physical} centroid of cell $(i,j)$.
This step is linearity preserving.
\begin{table}[H]
\centering
\caption{\textbf{Linear Redistribution Formulas}}
\label{eqn:back2}
\begin{tabular}{c c}
\toprule
method  &  redistribution formula \\[.08in]
SRD & $\displaystyle U_{i,j}^{n+1} = \frac{1}{N_{i,j}}  \sum_{(r,s)\in W_{i,j}} \,
{q}^{srd}_{r,s}(x_{i,j},y_{i,j}) \qquad\quad $ \\[.17in]
DRD & $\displaystyle U_{i,j}^{n+1} = U_{i,j}^n -  \frac{1}{N_{i,j}}  \sum_{(r,s)\in W_{i,j}} \,
q^{drd}_{r,s}(x_{i,j},y_{i,j}) $\\[.08in]
\bottomrule
\end{tabular}
\end{table}
\end{itemize}

\small
\bibliography{references}
\bibliographystyle{plain}

\end{document}